\documentclass[a4paper, 11pt]{amsart}
\usepackage{amsmath,amssymb,amscd}
\usepackage{amsfonts}
\usepackage[english]{babel}
\usepackage[leqno]{amsmath}
\usepackage{amssymb,amsthm}
\usepackage{mathrsfs} 
\usepackage{amscd}
\usepackage{enumerate}
\usepackage{epsfig}
\usepackage{relsize}
\usepackage{layout}
\usepackage{fullpage}
\usepackage[usenames,dvipsnames]{xcolor}
\usepackage[backref=page]{hyperref}
\usepackage{tikz}
\usepackage{tikz-cd}
\usetikzlibrary{matrix,arrows}

\hypersetup{
 colorlinks,
 citecolor=Green,
 linkcolor=Red,
 urlcolor=Blue}
\usepackage[matrix,arrow,tips,curve]{xy}
\input{xy}
\xyoption{all}
\definecolor{darkgreen}{rgb}{0.0, 0.7, 0.0}
\definecolor{purple}{rgb}{0.5, 0.0, 0.5}
\definecolor{red}{rgb}{0.8, 0.2, 0.0}

\newtheorem*{thm*}{Theorem}
\newtheorem{thm}{Theorem}[section]
\newtheorem{lemma}[thm]{Lemma}

\newtheorem{prop}[thm]{Proposition}

\numberwithin{equation}{section}
\theoremstyle{definition}
\newtheorem{defi}[thm]{Definition}

\newtheorem*{question*}{Question}
\newtheorem*{ack*}{Acknowledgments}
\newtheorem{notation}[thm]{Notation}

\theoremstyle{remark}
\newtheorem{remark}[thm]{Remark}

\newtheorem{example}[thm]{Example}

\DeclareMathOperator{\Supp}{Supp}
\DeclareMathOperator{\Ker}{Ker}
\DeclareMathOperator{\Sym}{Sym}
\DeclareMathOperator{\Ima}{Im}
\DeclareMathOperator{\Gr}{Gr}
\DeclareMathOperator{\Coker}{Coker}

\title{Partially ample line bundles and base loci}

\author{Marco D'Ambra}

\address{\hskip -.43cm Dipartimento di Matematica, Universit\'a di Roma, Tor Vergata, V.le della Ricerca Scientifica, I-00133 Roma, Italy, e-mail {\tt marcodambra91@gmail.com}}

\begin{document}


\begin{abstract}
We generalize some results of A.J. Sommese, B. Totaro and M.V. Brown, providing some geometric interpretations of the notion of partial ampleness.
\end{abstract}

\maketitle

\section*{Introduction} \label{intro}

The most important notion of positivity in algebraic geometry is undoubtedly the ampleness.
By the well-known Cartan-Grothendieck-Serre's result, the ampleness of a line bundle is characterized in terms of the vanishing of the higher cohomology groups. By weakening this condition we obtain a notion of partial ampleness, that intuitively measures how much a line bundle is far from being ample and that shares many important properties with traditional ampleness.

Let $X$ be a projective scheme of dimension $n$, let $q=0,\dots ,n-1$ be an integer and let $L$ be a line bundle on $X$. Then $L$ is called $q$-ample if for all coherent sheaves $\mathcal{F}$ we have that $H^i(X,L^m\otimes \mathcal{F})=0$ for all $m\gg 0$ and $i >q$.
Observe that for $q=0$ we recover the notion of ampleness.

Partial ampleness has been studied, among others, by Sommese, Totaro, K\"uronya and Ottem (c.f. \cite{Sommese,Totaro,KuronyaACF,KuronyaPS,Ottem})

Even if the definition is purely cohomological, there exist in literature some interesting results that help us to interpret geometrically the notion of partial ampleness. One of the most important ones is due to Sommese (c.f. \cite[Theorem 1.7]{Sommese}, see Theorem \ref{semiamplecase}), who proved that a semiample line bundle is $q$-ample if and only if the fibers of its Iitaka fibration have dimension at most $q$. We provide a generalization of this result by relaxing the hypotheses on $L$.
Assume that $k(X,L)\geq 0$ and let $m_0>0$ be an integer such that $\rm{Bs}|m_0L|$ is the stable base locus $\mathbf{B}(L)$. Moreover let $\pi:\hat{X}\rightarrow X$ be the blow-up of the base ideal of $|m_0L|$, with exceptional divisor $E$. We have the decomposition $|\pi^*(m_0L)|=|M|+E$, with $M$ base-point-free. 
Sommese's result characterizes the partial ampleness of $M$. Moreover, if $L$ is assumed to be semiample, then $m_0L=M$. It is then natural to compare the partial ampleness of $L$ and $M$. We have the following result:
\begin{thm*}[Theorem \ref{LvsM}]
In the previous setting:
\begin{itemize}
	\item[(i)] If $\dim \mathbf{B}(L)\leq q$, then $M$ $q$-ample implies $L$ $q$-ample.
	\item[(ii)] If $\dim \mathbf{B}(L)\leq q-1$, then $M$ is $q$-ample if and only if $L$ is $q$-ample.
\end{itemize}
\end{thm*}
We also provide examples that show that the hypotheses on $\dim \mathbf{B}(L)$ are sharp.\\

Another beautiful geometric interpretation of partial ampleness is due to Totaro, who proved that a line bundle $L$ is $(n-1)$-ample if and only if $-L$ is not pseudoeffective (c.f. \cite[Theorem 9.1]{Totaro}, see Theorem \ref{Totaro}). More recently Brown proved that if $L$ is big, then $L$ is $(n-2)$-ample if and only if for every subvariety $Z$ of dimension $n-1$ we have that $-L_{|Z}$ is not pseudoeffective (\cite[Corollary 1.2]{Brown}, see Theorem \ref{Brown}). We generalize these two results with the following theorem, which translates the $q$-ampleness of a line bundle $L$ with augmented base locus $\mathbf{B}_+(L)$ of dimension at most $q+1$ in terms of the geometry of its restriction to the subvarieties of $X$:
\begin{thm*}[Theorem \ref{Brownies}]
Let $X$ be a pure dimensional scheme and let $L$ be a line bundle on $X$. If $\dim \mathbf{B}_+(L)\leq q+1$, then the following conditions are equivalent:
\begin{itemize}
	\item[(i)] $L$ is $q$-ample.
	\item[(ii)] For all subvarieties $Z\subset X$ of dimension $q+1$ we have that $-L_{|Z}$ is not pseudoeffective.
\end{itemize}
\end{thm*}
Observe that if $q=n-1$ the hypothesis on $\dim \mathbf{B}_+(L)$ is empty and we recover Totaro's result, while if $q=n-2$ the hypothesis on $\dim \mathbf{B}_+(L)$ is equivalent to the bigness of $L$ and we generalize Brown's result.\\

The pull-back of an ample line bundle under a blow-up is never ample. Similarly the pull-back of a $q$-ample line bundle in never $q$-ample. However, the following theorem shows that we can expect some partial ampleness.
\begin{thm*}[Theorem \ref{cblowup}]
Let $X$ be a projective normal variety, let $Z$ be an irreducible l.c.i. subscheme of codimension $e\geq 1$, let $\pi: \hat{X} \rightarrow X$ be the blow-up of $X$ along $Z$ and let $L$ be a line bundle on $X$.
If $L$ is $q$-ample, then $\pi^*L$ is $(q+e-1)$-ample.
\end{thm*}

The paper is organized as follows: in section \ref{def} we give a quick review of the theory of partial ampleness, in section \ref{sa} we prove Theorem \ref{LvsM}, in section \ref{sub} we prove Theorem \ref{Brownies} and in section \ref{last} we prove Theorem \ref{cblowup}.

\begin{ack*}
This work is part of my PhD thesis. I would like to thank my advisor, Angelo Lopez, for his invaluable guidance.
\end{ack*}

\section{Definitions and first properties} \label{def}

\begin{notation}
Unless otherwise specified $X$ will be a projective noetherian scheme of finite type and of dimension $n$ over the complex number field and $q=0,\dots ,n-1$ will be an integer. With the term (sub-)variety we mean a reduced and irreducible (sub-)scheme.
\end{notation}
 
Let $L$ be a line bundle on $X$.
The stable base locus of $L$ is the Zariski closed subset
$$\mathbf{B}(L) =\bigcap_{m\in \mathbb{N}}\rm{Bs}|mL|.$$
The augmented base locus of $L$ is
$$\mathbf{B}_+(L)=\bigcap_{m\in \mathbb{N}}\mathbf{B}(mL-A)$$
where $A$ is an ample line bundle.\\
For more information on stable base loci and augmented base loci we refer to \cite{ELMNPBL}.

\begin{defi}
$L$ is $q$-ample if for all coherent sheaves $\mathcal{F}$ on $X$ there exists an integer $m_{L,\mathcal{F}}>0$ such that
$$H^i(X,L^m\otimes \mathcal{F})=0$$
for all $m> m_{L,\mathcal{F}}$ and $i >q$.
\end{defi}

By the definiton it follows immediately that, if $L$ is $q$-ample, then it is $(q+1)$-ample.
Moreover if $L$ is $q$-ample, then $mL$ is $q$-ample for all positive integer $m>0$.
Observe finally that for $q=0$, by the Cartan-Grothendieck-Serre's result, we recover the notion of ampleness.

The following theorem is the most important characterization of $q$-ample line bundles:

\begin{thm}[{\cite[Theorem 7.1]{Totaro}}]\label{cequivalentconditions}
Let $A$ be an ample divisor on $X$. Then there exists a constant $C_{X,A,q}>0$ such that
for all line bundles $L$ on $X$ the following conditions are equivalent:
\begin{itemize}
	\item[(i)] $L$ is $q$-ample.
	\item[(ii)] There exists an integer $m_{L,A}>0$ such that
	$$H^i(X,m_{L,A}L-rA)=0$$
	for all $1\leq r\leq C_{X,A,q}$ and $i>q$.
\end{itemize}
\end{thm}

The next proposition summarizes the most important properites that the notion of partial ampleness shares with the usual one:

\begin{prop}[{\cite[ Corollary 7.2]{Totaro}, \cite[Proposition 2.3]{Ottem}}] \label{credirr}
Let $L$ be a line bundle on $X$. Then
\begin{itemize}
	\item[(i)] $L$ is $q$-ample if and only if $L_{red}$ is so.
	\item[(ii)] $L$ is $q$-ample if and only if $L_{|X_i}$ is $q$-ample on every irreducible component $X_i$ of $X$.
	\item[(iii)] Let $f:Y\rightarrow X$ be a finite morphism of projective schemes.
Then if $L$ is $q$-ample we have that $f^*L$ is also $q$-ample. If $f$ is surjective, then the converse also hold.
In particular, if $Z\subset X$ is a subscheme and $L$ is $q$-ample, then $L_{|Z}$ is also $q$-ample.
\end{itemize}
\end{prop}

The following proposition relates the $q$-ampleness of a line bundle and the dimension of its augmented base locus:

\begin{prop}[{\cite[Theorem B, Corollary 2.6]{KuronyaPS}}] \label{KuronyaPS-B}
Let $L$ be a line bundle on $X$. Then for all coherent sheaves $\mathcal{F}$ on $X$ there exists an integer $m_{L,\mathcal{F}}>0$ such that
$$H^i(X,L^m\otimes D \otimes \mathcal{F})=0 $$
for all nef divisors $D$, $m>m_{L,\mathcal{F}}$ and $i>\dim  \mathbf{B}_+(L)$ .
In particular, if $\dim \mathbf{B}_+(L)\leq q$, then $L$ is $q$-ample.
\end{prop}

The next example shows that, even in the semiample case, it is not true that a $q$-ample line bundle has the augmented base locus of dimension at most $q$, disproving a result of Choi (c.f. \cite[Theorem 1.1]{Choi}):

\begin{example}\label{antichoi}
Let $\pi:X\rightarrow Y$ be the blow-up of a smooth threefold $Y$ along a smooth curve $C$ with exceptional divisor $E$ and let $H$ be a very ample divisor on $Y$.
Since $H$ is ample and $\pi$ is birational, we have by \cite[Proposition 2.3]{BBP} that
$$ \mathbf{B}_+(\pi^* H) = \pi^{-1}( \mathbf{B}_+(H))\cup E =E.$$
Hence $\dim \mathbf{B}_+(\pi^*H)=2$.
Moreover, since $H$ is very ample, we have that $\pi^*H$ is semiample, $\mathbb{P}H^0(X,\pi^*H)=\mathbb{P}H^0(Y,H)$ and the morphism $\phi_{|\pi^* H|}$ factorizes through the blow-up $\pi$. Observe now that $\phi_{|H|}$ is an embedding and that the fiber of $\pi$ over a point in $Y$ can be a point or a $\mathbb{P}^1$. Hence the dimension of the fibers of the $\phi_{|\pi^*H|}$ is at most $1$ and so by \cite[Theorem 1.7]{Sommese} (see Theorem \ref{semiamplecase}) we have that $\pi^*H$ is $1$-ample.
\end{example}

We conclude this section by recalling these two theorems, which give a characterization of partial ampleness in the cases $q=n-1$ and $q=n-2$:

\begin{thm}[{\cite[Theorem 9.1]{Totaro}}]  \label{Totaro}
Let $L$ be a line bundle on $X$. Then the following conditions are equivalent:
\begin{itemize}
	\item[(i)] $L$ is $(n-1)$-ample.
	\item[(ii)] $-L$ is not pseudoeffective.
\end{itemize}
\end{thm}

\begin{thm}[{\cite[Corollary 1.2]{Brown}}] \label{Brown}
Let $X$ be a smooth variety and let $L$ be a big line bundle on $X$. Then the following conditions are equivalent:
\begin{itemize}
	\item[(i)] $L$ is $(n-2)$-ample.
	\item[(ii)] For all subvariety $Z\subset X$ of dimension $n-1$ we have that $-L_{|Z}$ is not pseudoeffective.
\end{itemize}
\end{thm}

\section{Partial ampleness and semiample fibration} \label{sa}

The following theorem characterizes the $q$-ampleness of a semiample line bundle in terms of the dimension of the fibers of the semiample fibration:

\begin{thm}[{\cite[Theorem 1.7]{Sommese}, \cite[Theorem 2.44]{GrebKuronya}}] \label{semiamplecase}
Let $L$ be a semiample line bundle on $X$ and let $m_0>0$ be a positive integer such that $m_0L$ is base-point-free. Then the following conditions are equivalent:
\begin{itemize}
	\item[(i)] $L$ is $q$-ample.
	\item[(ii)] The dimension of the fibers of the morphism $\phi_{|m_0L|}$ is at most $q$.
	\end{itemize}
\end{thm}




Inspired by Theorem \ref{semiamplecase} we have proved this result:

\begin{thm}\label{LvsM}
Let $L$ be a line bundle on $X$ with $k(X,L)\geq 0$ and let $m_0>0$ be an integer such that $\mathbf{B}(L)=\rm{Bs}|m_0L|$.
Moreover let $\pi:\hat{X}\rightarrow X$ be the blow-up of the base ideal of $|m_0L|$, with exceptional divisor $E$ 
and denote $M=\pi^*(m_0L)-E$. Then:
\begin{itemize}
	\item[(i)] If $\dim \mathbf{B}(L)\leq q$, then $M$ $q$-ample implies $L$ $q$-ample.
	\item[(ii)] If $\dim \mathbf{B}(L)\leq q-1$, then $M$ is $q$-ample if and only if $L$ is $q$-ample.
\end{itemize}
\end{thm}

\proof
Let $A$ be an ample line bundle on $X$ such that $\hat{A}=\pi^*A-E$ is ample on $\hat{X}$ and let $C\geq 1$ be a constant.
Moreover take $m,t\geq 1$. We want to compare the cohomology groups $H^i(X,mm_0L-tA)$ and $H^i(\hat{X},mM-t\hat{A})$ with the Leray spectral sequence
$$E_1^{p,r}=H^p(X,R^{r-p}\pi_*(mM-t\hat{A}))\Rightarrow H^r(\hat{X},mM-t\hat{A})$$
where for spectral sequences we are using the notations of \cite{HS}.
By projection formula we get that
$$E_1^{p,r}= H^p(X,R^{r-p}\pi_*(\pi^*(mm_0L-tA)+(t-m)E))=$$
$$=H^p(X,\mathcal{O}_X(mm_0L-tA)\otimes R^{r-p}\pi_*\mathcal{O}_{\hat{X}}((t-m)E)).$$
Since $-E$ is $\pi$-ample, it follows by \cite[Theorem 1.7.6]{LazPAG}  that for all $1\leq t\leq C$ there exists an integer $m_t>0$ such that 
$$R^i\pi_*\mathcal{O}_{\hat{X}}( (t-m)E)=0$$
for all $i>0$ and $m\geq m_t$.
Taking $M_C=\max\{m_1,\dots,m_C,C+1\}$ we have that
\begin{equation}\label{eq:nome}
E_1^{p,r} =
\begin{cases}
	0 & \text{if } r-p\neq 0  \\
	H^p(X,\mathcal{O}_X(mm_0L-tA)\otimes \pi_*\mathcal{O}_{\hat{X}}( (t-m)E) ) & \text{if } r-p = 0
	\end{cases}
\end{equation}
for all $m\geq M_C$ and $1\leq t\leq C$.\\
If $r-p\neq 0$ we have that $E_\infty^{p,r}=E_1^{p,r}=0$. If $r-p=0$ consider the maps
$$E_l^{r+l+1,r+1}\rightarrow E_l^{r,r}\rightarrow E_l^{r-l-1,r-1}$$
with $l\geq 1$.
Since $E_l^{r+l+1,r+1}=E_l^{r-l-1,r-1}=0$ for all $l\geq 1$ we get that $E_{\infty}^{r,r}=E_1^{r,r}$.\\
Consider now an integer $i\geq 0$ and the filtration
$$H^i(\hat{X},mM-t\hat{A})=F^0\supset F^1\supset \cdots \supset F^k\supset F^{k+1}=0$$
with $F^p/F^{p+1}=\Gr^p(H^i(\hat{X},\pi^*(mm_0L)-t\hat{A}))=E_\infty^{p,i}$ for all $p\geq 0$.\\
Since $F^p/F^{p+1}=0$ for all $p\neq i$ the filtration is
$$H^i(\hat{X},mM-t\hat{A})=F^0= F^i\supset F^{i+1}=0.$$
By the exact sequence
$$0\rightarrow F^{i+1} \rightarrow F^i\rightarrow F^i/F^{i+1}\rightarrow 0 $$
and by (\ref{eq:nome}) we get that
\begin{equation}\label{eq:nome2}
H^i(\hat{X},mM-t\hat{A}) \cong H^i(X,\mathcal{O}_X(mm_0L-tA)\otimes \pi_*\mathcal{O}_{\hat{X}}((t-m)E))
\end{equation}
for all $i\geq 0$, $m\geq M_C$ and $1\leq t\leq C$.\\
Now we prove $(i)$. Assume that $M$ is $q$-ample and that $\dim \mathbf{B}(L)\leq q$. Let $C=C_{X,A,q}$ be the constant of Theorem \ref{cequivalentconditions}.
To show that $L$ is $q$-ample we prove that $m_0L$ is $q$-ample by finding an integer $m_L>0$ such that
$$H^i(X,m_Lm_0L-tA)=0$$
for all $i>q$ and $1\leq t\leq C$.\\
To see this fix $i>q$, $m\geq M_C$ and $1\leq t \leq C$.
Since $m>t$, then $\pi_*\mathcal{O}_{\hat{X}}( (t-m)E)$ is an ideal sheaf supported on $\mathbf{B}(L)$. Set now $Z= Z(\pi_*\mathcal{O}_{\hat{X}}( (t-m)E))$ and consider the exact sequence
$$0\rightarrow \mathcal{O}_X(mm_0L-tA)\otimes \pi_*\mathcal{O}_{\hat{X}}( (t-m)E)\rightarrow \mathcal{O}_X(mm_0L-tA)\rightarrow \mathcal{O}_{Z}(mm_0L-tA)\rightarrow 0 $$
By (\ref{eq:nome2}) we get that
$$h^i(X,mm_0L-tA)\leq h^i(\hat{X}, mM-t\hat{A})+h^i(Z,(mm_0L-tA)_{|Z}). $$
Since $\dim \mathbf{B}(L)\leq q$ the second term on the right is zero and we have the inequality
$$h^i(X,mm_0L-tA)\leq h^i(\hat{X}, mM-t\hat{A}).$$
Since $M$ is $q$-ample we have that for all $1\leq t \leq C$ there exists an integer $\hat{m}_t>0$ such that
$$H^i(X,mM-t\hat{A})=0$$
for all $i>q$ and $m\geq\hat{m}_t$.
Taking $m_L=\max\{\hat{m}_1,\dots ,\hat{m}_C, M_C\}$ we have that 
$$H^i(X,m_Lm_0L-tA)=0$$
for all $i>q$ and $1\leq t\leq C$. This proves $(i)$.\\
The proof of $(ii)$ is similar. Assume that $L$ is $q$-ample and that $\dim \mathbf{B}(L)\leq q-1$. Let $C=C_{\hat{X}, \hat{A}, q}$ be the constant of Theorem \ref{cequivalentconditions}. To show that $M$ is $q$-ample we find an integer $m_M>0$ such that
$$H^i(\hat{X},m_MM-t\hat{A})=0$$
for all $i>q$ and $1\leq t\leq C$.\\
Take now $i>q$, $m\geq M_C$ and $1\leq t \leq C$.
For all integers $s\geq 1$ we have an exact sequence
$$ 0\rightarrow \mathcal{O}_{\hat{X}} ((t+s-1-m)E) \rightarrow \mathcal{O}_{\hat{X}} ((t+s-m)E)\rightarrow \mathcal{O}_E((t+s-m)E)\rightarrow 0 $$
Applying the functor $\pi_*$ we get the exact sequence
$$ 0\rightarrow \pi_* \mathcal{O}_{\hat{X}} ((t+s-1-m)E) \xrightarrow{\varphi_s} \pi_*\mathcal{O}_{\hat{X}} ((t+s-m)E)\rightarrow\mathcal{F}_s \rightarrow 0$$
where $\mathcal{F}_s =\Coker \varphi_s$ is a coherent sheaf supported on $\mathbf{B}(L)$. Tensoring by $\mathcal{O}_X(mm_0L-tA)$ we get the exact sequence
$$0\rightarrow \mathcal{O}_X(mm_0L-tA)\otimes \pi_*\mathcal{O}_{\hat{X}}((t+s-1-m)E) \rightarrow \mathcal{O}_X(mm_0L-tA) \otimes 
\pi_*\mathcal{O}_{\hat{X}}((t+s-m)E)$$
$$\rightarrow \mathcal{O}_{X}(mm_0L-tA) \otimes \mathcal{F}_s \rightarrow 0 $$
By (\ref{eq:nome2}) we have that
$$h^i(\hat{X},mM-t\hat{A}) = h^i(X,\mathcal{O}_X(mm_0L-tA) \otimes \pi_*\mathcal{O}_{\hat{X}}( (t-m)E))\leq  $$
$$\leq h^i(X,\mathcal{O}_X(mm_0L-tA))+\sum_{s=1}^{m-t} h^{i-1}(\mathbf{B}(L),\mathcal{O}_{X}(mm_0L-tA)\otimes \mathcal{F}_s).$$
Since $\dim\mathbf{B}(L)\leq q-1$, the second term on the right hand side is zero.
Thus we have the inequality 
$$h^i(\hat{X}, mM-t\hat{A})\leq h^i(X,mm_0L-tA).$$
Since $L$ is $q$-ample, for all $1\leq t\leq C$ there exists an integer $\hat{m}_t>0$ such that
$$H^i(X,mm_0L-tA)=0$$
for all $i>q$ and $m\geq\hat{m}_t$.
Taking $m_M=\max\{\hat{m}_1,\dots ,\hat{m}_C,M_C\}$ we conclude.
\endproof

The next example shows that the first assertion of Theorem \ref{LvsM} is sharp. Namely, if $M$ is $q$-ample and $\dim \mathbf{B}(L)=q+1$, it is not always true that $L$ is $q$-ample.

\begin{example}\label{exK3}
Let $X\subset \mathbb{P}^3$ be a smooth surface of degree $d\geq 4$ containing a line $C$ and let $H$ be the hyperplane section. Observe that $C=\mathbf{B}(C)$ and that $C^2=2-d$ and $H.C=1$.\\
Consider now the divisor $L=H+C$. We claim that $C = \mathbf{B}(L)=\mathbf{B}_+(L)$.
Indeed since $L.C=3-d$ we have that $C \subseteq \mathbf{B}(L)$. 
On the other hand $\mathbf{B}_+(L)\subseteq \mathbf{B}(L-H) =\mathbf{B}(C)=C$.
Hence $\dim\mathbf{B}_+(L)=1 $ and so by Proposition \ref{KuronyaPS-B} we conclude that $L$ is $1$-ample. However, since $L.C=3-d$, it is not ($0$-)ample.\\
Let now $\pi:\hat{X}\rightarrow X$ be the blow-up of $X$ along $C$. Since $\pi$ is an isomorphism we have the decomposition $L=M+F$, with $M=H$ and $F=C$. Hence $M$ is ($0$-)ample and $\dim \mathbf{B}(L)=1$ but $L$ is only $1$-ample.
\end{example}

The following example shows that the second assertion of Theorem \ref{LvsM} is sharp. Namely, if $L$ is $q$-ample and $\dim \mathbf{B}(L)\geq q$, it is not always true that $M$ is $q$-ample.

\begin{example}
Let $p: X \rightarrow \mathbb{P}^n$ be the blow-up of $\mathbb{P}^n$ along a projective subspace $\mathbb{P}^{n-q-1}$ with exceptional divisor $L$. Since $\mathcal{O}_{\mathbb{P}^n}(1)$ is an ample line bundle, then $\mathcal{O}_{\mathbb{P}^n}(1)^{\oplus (q+1)}$ is an ample vector bundle of rank $q+1$. It follows by \cite[Proposition 4.5]{Ottem} that $\mathbb{P}^{n-q-1}$ is ample, i.e. $L$ is $q$-ample. Moreover we have that $ \mathbf{B}(L)=L$.
Let now $\pi:\hat{X}\rightarrow X$ be the blow-up of $X$ along $L$ with exceptional divisor $E$. Since $\pi$ is an isomorphism we have the decomposition $L=M+F$, with $M=0$ and $F=L$.
Since $-M=0$ is pseudoeffective, Theorem \ref{Totaro} implies that $M$ is not $(n-1)$-ample. On the other hand $L$ is $q$-ample and $\dim \mathbf{B}(L)=n-1\geq q$.
\end{example}

\section{Partial ampleness on subschemes} \label{sub}

In this section we provide a generalization of Theorem \ref{Totaro} and of Theorem \ref{Brown}.

\begin{lemma}\label{pluto}
Let $L$ be a pseudoeffective $\mathbb{R}$-divisor on $X$. Then for all very general hyperplane sections $H$ the restriction $L_{|H}$ is pseudoeffective.
\end{lemma}
\proof
Since $L$ is pseudoeffective, there exists a sequence $\{ D_m \}_{m\geq 1}$ of effective $\mathbb{R}$-divisors such that $ [L]=\lim_{m\rightarrow \infty} [D_m]$ in $N^1(X)$.
Since $\bigcup_{m\geq 1} \Supp D_m$ is (at most) a countable union of subschemes of dimension $n-1$, then we can take an hyperplane section $H$ such that 
for all $m\geq 1$ and $G\in \Supp (D_m)$ we have that $H\neq G$.
It follows that $[D_{m|H}]$ is effective for all $m\ge 1$ and so $[L_{|H}]=\lim_{m\rightarrow \infty} [D_{m|H}]$ is pseudoeffective.
\endproof

\begin{prop}\label{porcone}
Let $L$ be a line bundle on $X$.
The following conditions are equivalent:
\begin{itemize}
	\item[(i)] For all subvarieties $Z\subset X$ of dimension $q+1$ we have that $L_{|Z}$ is $q$-ample.	
	\item[(ii)] For all subvarieties $Z\subset X$ of dimension $q+1$ we have that $-L_{|Z}$ is not pseudoeffective.
	\item[(iii)] For all subvarieties $Z\subset X$ of dimension $>q$ we have that $-L_{|Z}$ is not pseudoeffective.
\end{itemize}
\end{prop}
\proof
To obtain $(i)\Leftrightarrow (ii)$ we observe that, since the dimension of $Z$ is $q+1$, we can apply Theorem \ref{Totaro}.
To get $(ii)\Rightarrow (iii)$ suppose that there exists a subvariety $Z$ of dimension $>q+1$ such that
 $-L_{|Z}$ is pseudoeffective. Taking a very general hyperplane section $H$ on $Z$ we have by Lemma \ref{pluto} that $-(L_{|Z})_{|H}$ is also pseudoeffective. If the dimension of $Z$ is $q+2$ we get a contradiction. Otherwise we can iterate the procedure until we find a subvariety $W$ of dimension $q+1$ such that $-L_{|W}$ is pseudoeffective.
Finally, the inclusion $(iii)\Rightarrow (ii)$ is obvious.
\endproof

The rest of the section is devoted to the proof that, under the addictional hypothesis that 
$\dim\mathbf{B}_+(L)\leq q+1$, the conditions of the previous proposition are also equivalent to the $q$-ampleness of $L$.\\

We start recalling the following result:

\begin{lemma}[{\cite[Lemma 2.2]{Brown}}] \label{Brown-2.2}
Let $L$ be a line bundle and let $A$ be an ample divisor on $X$.
If $L$ is $q$-ample, then for all coherent sheaves $\mathcal{F}$ on $X$ there exist two integers $m_{L,A,\mathcal{F}}, m'_{L,A,\mathcal{F}}>0$ such that
$$H^i(X,L^m\otimes A^{m'}\otimes \mathcal{F})=0$$
if $m>m_{L,A,\mathcal{F}}$ or $m'>m'_{L,A,\mathcal{F}}$ and $i>q$.
\end{lemma}

The next theorem of Brown is crucial in the proof of our characterization.
For completeness, we provide a quite different version of the proof:

\begin{thm}[{\cite[Theorem 2.1]{Brown}}]\label{Brown-2.1+}
Let $L$ and $L'$ be two line bundles on $X$ and let $a,b>0$ be two integers such that $h^0(X,L')>0$ and $A:=aL-bL'$ is ample. Moreover let $s\in H^0(X,L')$ be a non-zero section of $L'$ and let $Z(s)$ be the associated subscheme.
If $L$ is not $q$-ample, then $L_{|Z(s)}$ is not $q$-ample.
\end{thm}
\proof
First of all we make some reductions.\\
We may assume that $a=1$.
Indeed, since a line bundle is $q$-ample if and only if every positive multiple is $q$-ample, we have that $aL$ is not $q$-ample.
Moreover if we show that $aL_{|Z(s)}$ is not $q$-ample, then $L_{|Z(s)}$ is also not $q$-ample.\\
We may assume that $b=1$. To see this take the non-zero section $s^{\otimes b}\in H^0(X,bL')$.
We have that $Z(s^{\otimes b})_{red}=Z(s)_{red}$.
Thus if we show that $L_{|Z(s^{\otimes b})}$ is not $q$-ample, then, by Proposition \ref{credirr} part $(i)$, $L_{|Z(s)}$ is also not $q$-ample. From now on we assume $a=b=1$ and so $L'=L-A$.\\
We may assume that $L$ is $(q+1)$-ample. 
To see this consider the integer
$$q_0=\max\{l\geq 0 : L \text{ is not $l$-ample} \}$$
and observe that $q\leq q_0< n$.
By definition $L$ is $(q_0+1)$-ample but not $q_0$-ample.
If we show that $L_{|Z(s)}$ is not $q_0$-ample, then it is also not $q$-ample.
Thus we can replace $q$ with $q_0$.\\
Consider now the exact sequence
$$0 \rightarrow \mathcal{F} \rightarrow -L'\xrightarrow{\varphi} \mathcal{O}_X \rightarrow O_{Z(s)} \rightarrow 0$$
for some sheaf $\mathcal{F}$, that is non-zero if $s$ is zero on some irreducible component of $X$. Let $\mathcal{G}=\Ima\varphi$.
Since $X$ is noetherian by \cite[Proposition II.5.7]{Hartshorne} we obtain that  $\mathcal{F}$ and $\mathcal{G}$ are coherent sheaves.
We have two short exact sequences
$$0 \rightarrow \mathcal{F} \rightarrow -L'\rightarrow \mathcal{G}\rightarrow 0,\ \ 0 \rightarrow \mathcal{G} \rightarrow \mathcal{O}_X\rightarrow \mathcal{O}_{Z(s)}\rightarrow 0.$$
We assume that $L_{|Z(s)}$ is $q$-ample and we reach a contradiction.\\
Denote by $C=C_{X,A,q}$ the constant of Theorem \ref{cequivalentconditions}.
Since $L$ is $(q+1)$-ample by Lemma \ref{Brown-2.2} there exists an integer $m_1>0$ such that
\begin{equation}\label{eq:ciao11}
H^i(X, \mathcal{F}(-CA)(mL+bA))=0
\end{equation}
for all $i>q+1$, $m> m_1$ and $b\geq 0$.
Moreover, since $L_{|Z(s)}$ is $q$-ample, again by Lemma \ref{Brown-2.2} there exists an integer $m_2>0$ such that
\begin{equation}\label{eq:ciao22}
H^i(Z(s), \mathcal{O}_{Z(s)}(-CA_{|Z(s)})((mL+bA)_{|Z(s)}))=0
\end{equation}
for all $i>q$, $m> m_2$ and  $b\geq 0$. Take $M_1=\max\{m_1,m_2\}$.\\
Since $A$ is ample there exists an integer $m_3>0$ such that
\begin{equation}\label{eq:ciao33}
H^i(X,(m-r)A+kL')=0
\end{equation}
for all $i>0$, $m> m_3$, $1\leq r\leq C$ and $0\leq k \leq M_1$.
Moreover, since $L$ is $(q+1)$-ample, by Theorem \ref{cequivalentconditions} there exists an integer $m_4>0$ such that
\begin{equation}\label{eq:ciao66}
H^i(X,mL-rA)=0
\end{equation}
for all $i>q+1$, $m> m_4$ and $1\leq r \leq C$. Take $M_2=\max\{m_3,m_4  \}$.\\
Since $L$ is not $q$-ample there exist an $i_0>q$, an $m_0> M_2$ and an integer $r_0$ such that $1\leq r_0 \leq C$ and
$$H^{i_0}(X, m_0L-r_0A)\neq 0.$$
Indeed otherwise we would have that
$$H^i(X,mL-rA)=0$$
for all $i>q$, $m> M_2$ and $1\leq r\leq C$ and hence that $L$ is $q$-ample by Theorem \ref{cequivalentconditions}.
Since $m_0> M_2\geq m_4$ by (\ref{eq:ciao66}) we have that
$$H^i(X, m_0L-r_0A)=0$$
for all $i>q+1$. Hence $i_0=q+1$ and 
\begin{equation}\label{eq:ciao44}
H^{q+1}(X,m_0L-r_0A)\neq 0.
\end{equation}
Consider now the set
$$P=\{l\in \mathbb{N}\ :\ H^{q+1}(X,(m_0-r_0)A+kL')=0\ \forall\   0\leq k\leq l\}.$$
Since $m_0> M_2\geq m_3$ by (\ref{eq:ciao33}) we get that $l \in P$ for all $0\leq l\leq M_1$.
Moreover by (\ref{eq:ciao44}) we have that $l \not\in P$ for all $l\geq m_0$. Hence there is a well-defined $k_0=\max P+1$ such that
\begin{equation}\label{eq:ciao55}
H^{q+1}(X,(m_0-r_0)A+(k_0-1)L' )=0,\ \ H^{q+1}(X,(m_0-r_0)A+k_0L' )\neq 0.
\end{equation}
Consider the exact sequence
$$0 \rightarrow \mathcal{F}((m_0-r_0)A+k_0L') \rightarrow (m_0-r_0)A+(k_0-1)L'\rightarrow \mathcal{G}((m_0-r_0)A+k_0L')\rightarrow 0$$
Since $k_0> M_1\geq m_1$ and $m_0-r_0-k_0+C\geq 0$ by (\ref{eq:ciao11}) we have that
$$H^{q+2}(X,\mathcal{F}((m_0-r_0)A+k_0L'))=H^{q+2}(X,\mathcal{F}((m_0-r_0-k_0)A+k_0L))=$$
$$=H^{q+2}(X,\mathcal{F}(-CA)((m_0-r_0-k_0+C)A+k_0L))=0.$$
It follows by (\ref{eq:ciao55})  that
\begin{equation}\label{eq:ciao77}
H^{q+1}(\mathcal{G}((m_0-r_0)A+k_0L'))=0.
\end{equation}
Consider now the exact sequence
$$0 \rightarrow \mathcal{G}((m_0-r_0)A+k_0L') \rightarrow (m_0-r_0)A+k_0L' \rightarrow ((m_0-r_0)A+k_0L')_{|Z(s)}\rightarrow 0$$
Since $k_0> M_1\geq m_2$ and $m_0-r_0-k_0+C\geq 0$ by (\ref{eq:ciao22}) we have that
$$H^{q+1}(Z(s),((m_0-r_0)A+k_0L')_{|Z(s)})=H^{q+1}(Z(s),((m_0-r_0-k_0)A+k_0L)_{|Z(s)})=$$
$$=H^{q+1}(Z(s),\mathcal{O}_{Z(s)}(-CA_{|Z(s)})(((m_0-r_0-k_0+C)A+k_0L)_{|Z(s)})=0.$$
It follows that
$$H^{q+1}((m_0-r_0)A+k_0L')=0.$$
This contradicts (\ref{eq:ciao55}) and so $L_{|Z(s)}$ is not $q$-ample.
\endproof

The following lemma, which will not be used later, shows that the hypothesis of the previous theorem is in fact equivalent to the bigness of $L$.

\begin{lemma}\label{paperoga}
Let $L$ be a line bundle on $X$. Then the following conditions are equivalent:
\begin{itemize}
	\item[(i)] $L$ is big, i.e. there exist an ample $\mathbb{Q}$-divisor $A$ and an effective $\mathbb{Q}$-divisor $E$ such that $L\sim_{\mathbb{Q}} A+E$.
	\item[(ii)] There exist a line bundle $L'$ and two integers $a,b>0$ such that $h^0(X,L')>0$ and $aL-bL'$ is ample.
\end{itemize}
\end{lemma}
\proof
If $L$ is big, then there exist an ample $\mathbb{Q}$-divisor $A$ and an effective $\mathbb{Q}$-divisor $E$ such that $L\sim_{\mathbb{Q}} A+E$. Take $k>0$ such that $kL\sim_{\mathbb{Z}} kA+kE$, with $kA$ and $kE$ integer divisors and set $L'=kL-kA$. We get that $h^0(X,L')=h^0(X,kE)>0$. Taking $a=k$ and $b=1$ we get that $aL-bL'=kA$ is ample. This proves $(i) \Rightarrow (ii)$.\\
To prove $(ii) \Rightarrow (i)$ set $A=aL-bL' $ and observe that $L=\frac{1}{a} A + \frac{b}{a} L'$ is big because is a sum of an ample $\mathbb{Q}$-divisor and of an effective $\mathbb{Q}$-divisor.
\endproof

\begin{prop}\label{mylittleponyci}
Let $X$ be a projective pure dimensional scheme and let $L$ be a line bundle on $X$ such that $\dim\mathbf{B}_+(L)\leq q+1$.
If $L$ is not $q$-ample, then there exists a complete intersection subscheme $Z\subset X$ of pure dimension $q+1$ such that $L_{|Z}$ is not $q$-ample
\end{prop}
\proof
If $q=n-1$ take $Z=X$. If $0\leq q\leq n-2$ consider an ample divisor $A$ on $X$ and an integer $m_0>0$ such that
$$\mathbf{B}_+(L)=\rm{Bs}|m_0L-A|.$$
Denote $L'=m_0L-A$. For all $s\geq 0$, $1\leq \alpha \leq s$ and $s_\alpha \in H^0(X,L')$ denote $E=Z(s_\alpha)$ and set
$$
W_s=
\begin{cases}
	X & \text{if } s=0\\
	E_1\cap \cdots \cap E_s & \text{if } s\geq 1
\end{cases}.
$$
Moreover denote (adopting the notation of \cite{ELMNPRV})
$$H^0(X|W_s,L')=\Ima \Big(H^0(X,L')\rightarrow H^0(W_s, L'_{|W_s})  \Big).$$
We prove by descending induction on $q$ that there exist sections $s_1,\dots ,s_{n-q-1}\in H^0(X,L')$ such that:
\begin{itemize}
\item[(i)] $E_\alpha$ is an effective Cartier divisor on $X$ for all $1\leq \alpha \leq n-q-1$.
\item[(ii)] $Z=W_{n-q-1}$ is an effective Cartier divisor on $W_{n-q-2}$ of pure dimension $q+1$.
\item[(iii)] $L_{|Z}$ is not $q$-ample.
\end{itemize}
If $q=n-2$, since $X$ has pure dimension $n$ and $\mathbf{B}_+(L)=\rm{Bs}|L'|$ has dimension at most $q+1$, then $X_i\not\subset \rm{Bs}|L'|$ for all irreducible component $X_i$ of $X$.
Hence there exist points $x_i\in X_i\setminus \rm{Bs}|L'|$ and sections $s_{1,i}\in H^0(X,L')$ such that $s_{1,i}(x_i)\neq 0$.
It follows that $H^0(X,\mathcal{I}_{\{x_i\} / X}(L'))$ is properly contained in $H^0(X,L')$. Since we are dealing with vector spaces and the irreducible components of $X$ are finitely many we have that
\begin{equation}\label{eq:cazzoA}
\bigcup_iH^0(X,\mathcal{I}_{\{x_i\} / X}(L')) \subsetneq H^0(X,L').
\end{equation}
Hence there exists a section $s_1\in H^0(X,L')$ such that $s_{1|X_i}\neq 0$ for all $i$ and so $E_1$ is an effective Cartier divisor on $X$. Since $A$ is ample and $L$ is not $q$-ample, it follows by Theorem \ref{Brown-2.1+} that $L_{|E_1}$ is also not $q$-ample.\\
If $0\leq q\leq n-3$ assume by induction that there exist sections $s_1,\dots ,s_{n-q-2}\in H^0(X,L')$ such that $E_\alpha$ is an effective Cartier divisor on $X$ for all $1\leq \alpha \leq n-q-2$, $W_{n-q-2}$ is an effective Cartier divisor on $W_{n-q-3}$ of pure dimension $q+2$ and $L_{|W_{n-q-2}}$ is not $q$-ample.
Observe that 
$$\mathbf{B}_+(L)\cap W_{n-q-2}=\rm{Bs}|L'|\cap W_{n-q-2}=\rm{Bs}|H^0(X|W_{n-q-2},L')|.$$
Since $W_{n-q-2}$ has pure dimension $q+2$ and $\mathbf{B}_+(L)=\rm{Bs}|L'|$ has dimension at most $q+1$, then $Y_j\not\subset \rm{Bs}|L'|$ for all irreducible component $Y_j$ of $W_{n-q-2}$. As before there exist points $y_j\in Y_j\setminus \rm{Bs}|L'|$ and sections $s_{n-q-1,j}\in H^0(X,L')$ such that $s_{n-q-1,j}(y_j)\neq 0$.
It follows that
\begin{equation}\label{eq:cazzoB}
\bigcup_jH^0(W_{n-q-2},\mathcal{I}_{\{y_j\} / W_{n-q-2}}(L')) \cap H^0(X|W_{n-q-2},L') \subsetneq  H^0(X|W_{n-q-2},L').
\end{equation}
By (\ref{eq:cazzoA}) and (\ref{eq:cazzoB}) there exists an
$$s_{n-q-1} \in H^0(X,L')\setminus \bigcup_iH^0(X,\mathcal{I}_{\{x_i\} / X}(L'))$$
such that
$$s_{n-q-1|W_{n-q-2}}\in H^0(X|W_{n-q-2},L') \setminus \bigcup_jH^0(W_{n-q-2},\mathcal{I}_{\{y_j\} / W_{n-q-2}}(L')).$$
Hence $E_{n-q-1}$ is an effective Cartier divisor on $X$ and $Z=W_{n-q-1} = Z(s_1)\cap \cdots \cap Z(s_{n-q-1}) = Z(s_{n-q-1|W_{n-q-2}})$ is an effective Cartier divisor on $W_{n-q-2}$.
Since $A_{|W_{n-q-2}}$ is ample and $L_{|W_{n-q-2}}$ is not $q$-ample, then, by Theorem \ref{Brown-2.1+}, $L_{|Z}$ is also not $q$-ample and we conclude.

\begin{thm} \label{Brownies}
Let $X$ be a projective pure dimensional scheme and let $L$ be a line bundle on $X$ such that $\dim\mathbf{B}_+(L)\leq q+1$.
Then the following conditions are equivalent:
\begin{itemize}
	\item[(i)] $L$ is $q$-ample.
	\item[(ii)] For all subvarieties $Z\subset X$ of dimension $q+1$ we have that $L_{|Z}$ is $q$-ample.	
	\item[(iii)] For all subvarieties $Z\subset X$ of dimension $q+1$ we have that $-L_{|Z}$ is not pseudoeffective.
	\item[(iv)] For all subvarieties $Z\subset X$ of dimension $>q$ we have that $-L_{|Z}$ is not pseudoeffective.
\end{itemize}
\end{thm}
\proof
By Proposition \ref{porcone} we have $(ii)\Leftrightarrow (iii) \Leftrightarrow (iv)$.
To get $(i)\Rightarrow (ii)$ note that, by Proposition \ref{credirr}, for all subvarieties $Z\subset X$ of dimension $q+1$ we have that $L_{|Z}$ is $q$-ample.
To prove $(ii)\Rightarrow (i)$ assume that $L$ is not $q$-ample. Then by Proposition \ref{mylittleponyci}
there exists a complete intersection subscheme $Z'\subset X$ of pure dimension $q+1$ such that $L_{|Z'}$ is not $q$-ample.
By Proposition \ref{credirr} part $(i)$ there exists an irreducible component $Z''$ of $Z'$ such that $L_{|Z''}$ is not $q$-ample. Taking $Z=(Z'')_{red}$ we have by Proposition \ref{credirr} part $(ii)$ that $L_{|Z}$ is not $q$-ample.
\endproof

\begin{remark}
In particular, if $L$ is a line bundle such that $\dim\mathbf{B}_+(L)\leq 1$,
it follows by Theorem \ref{Brownies} that $L$ is ample if and only if it is strictly nef (i.e. $L.C>0$ for all curves $C$ in $X$).
\end{remark}

\section{Partial ampleness via blow-ups} \label{last}

In this section we prove the following result on the $q$-ampleness of the pull-back of a line bundle under a blow-up:

\begin{thm}\label{cblowup}
Let $X$ be a projective normal variety, let $Z$ be an irreducible l.c.i. subscheme of codimension $e\geq 1$, let $\pi: \hat{X} \rightarrow X$ be the blow-up of $X$ along $Z$ and let $L$ be a line bundle on $X$.
If $L$ is $q$-ample, then $\pi^*L$ is $(q+e-1)$-ample.
\end{thm}
\proof
Consider the commutative diagram
$$
\begin{CD}
E @>i>> \hat{X} \\
@V{p}VV @VV{\pi}V \\
Z @>>k> X \\
\end{CD}
$$
If $e=1$ we have that $\pi$ is an isomorphism and $\pi^*L$ is $q$-ample, so we may assume $e\geq 2$.\\
Consider an ample divisor $A$ on $X$ such that $\hat{A}=\pi^*A-E$ is ample on $\hat{X}$.
Moreover let $C=C_{\hat{X},\hat{A},q+e-1}$ be the constant of Theorem \ref{cequivalentconditions}.\\
To show that $\pi^*L$ is $(q+e-1)$-ample we find a 
constant $m_{\pi^*L}>0$ such that
$$H^i(\hat{X},\pi^*(m_{\pi^*L}L)-t\hat{A})=0$$
for all $i>q+e-1$ and $1\leq t\leq C$.\\
Set $i>q+e-1$ and $1\leq t \leq C$.\\
We want to compare the cohomology groups
$H^i(X,mL-tA)$ and $H^i(\hat{X},\pi^*(mL)-t\hat{A})$ using the Leray spectral sequence
$$E_1^{p,r}=H^p(X,R^{r-p}\pi_*(\pi^*(mL)-t\hat{A}))\Rightarrow H^r(\hat{X},\pi^*(mL)-t\hat{A}) $$
where $m\geq 1$.
By projection formula
$$E_1^{p,r}= H^p(X,R^{r-p}\pi_*(\pi^*(mL-tA)+tE))=H^p(X,\mathcal{O}_X(mL-tA)\otimes R^{r-p}\pi_*\mathcal{O}_{\hat{X}}(tE)).$$
Since $X$ is normal and $\pi$ is birational by (the proof of) \cite[Corollary III.11.4]{Hartshorne} we get
\begin{equation}\label{eq:1}
\pi_*\mathcal{O}_{\hat{X}}=\mathcal{O}_X.
\end{equation}
Moreover, using the theorem of formal functions, it can be proved as in \cite[Proposition V.3.4]{Hartshorne} that
\begin{equation}\label{eq:2}
R^s\pi_*\mathcal{O}_{\hat{X}}=0\ \ \forall\ s\geq 1.
\end{equation}
Since $Z$ is l.c.i. irreducible of codimension $e$ the normal bundle $N_{Z/X}$ is a vector bundle of rank $e-1$. Hence the exceptional divisor $E$ can be identified
with the projective space bundle $\mathbb{P}(N_{Z/X}^*)$ and $\mathcal{O}_E(E) = \mathcal{O}_{\mathbb{P}(N_{Z/X}^*)}(-1)$. 
Take now $l\geq 1$ and consider the exact sequence
$$0\rightarrow \mathcal{O}_{\hat{X}}((l-1)E)\rightarrow \mathcal{O}_{\hat{X}}(lE)\rightarrow \mathcal{O}_E(lE)\rightarrow 0 $$
Applying the functor $\pi_*$ we have a long exact sequence
$$0\rightarrow \pi_*\mathcal{O}_{\hat{X}}((l-1)E)\rightarrow \pi_*\mathcal{O}_{\hat{X}}(lE)\rightarrow 	\pi_*\mathcal{O}_E(lE)\rightarrow R^1\pi_*\mathcal{O}_{\hat{X}}((l-1)E)\rightarrow \cdots $$
Using the formulas for the direct images of $\mathcal{O}_E(lE)$ (see \cite[Appendix A]{LazPAG}) it's easy to see that
\begin{equation}\label{eq:3}
\pi_*\mathcal{O}_{\hat{X}}(lE)=\mathcal{O}_X\ \ \forall \ l\geq 1,
\end{equation}
\begin{equation}\label{eq:4}
R^s\pi_*\mathcal{O}_{\hat{X}}(lE)=0\ \ \forall\ l\geq 1,\ \ 1\leq s\neq e-1.
\end{equation}
It follows by (\ref{eq:1}), (\ref{eq:2}), (\ref{eq:3}) and (\ref{eq:4}) that, for all $m\geq 1$,
$$E_1^{p,r} =
\begin{cases}
	0 & \text{if } (r-p\leq -1) \lor ((r-p\geq 1)\land(r-p\neq e-1)  )\\
	H^p(X,mL-tA) & \text{if } r-p = 0\\
	H^p(X,\mathcal{O}_X(mL-tA)\otimes R^{e-1}\pi_*\mathcal{O}_{\hat{X}}(tE)) & \text{if } r-p=e-1
\end{cases}.
$$
If $(r-p\leq -1) \lor ((r-p\geq 1)\land(r-p\neq e-1)  )$ we have that $E_\infty^{p,r}=E_1^{p,r}=0$.\\
If $r-p=0$ consider the maps
$$E_l^{r+l+1,r+1}\rightarrow E_l^{r,r}\rightarrow E_l^{r-l-1,r-1}$$
with $l\geq 1$.
We have that $E_l^{r+l+1,r+1}=0$ for all $l\geq 1$ and that $E_l^{r-l-1,r-1}=0$ for all $1\leq l \neq e-1$.
Thus $E_{e-1}^{r,r}=E_1^{r,r}$ and $E_{\infty}^{r,r}=E_e^{r,r}$.
Looking at the maps $E_{e-1}^{r+e,r+1}\rightarrow E_{e-1}^{r,r}\rightarrow E_{e-1}^{r-e,r-1}$
and observing that $E_{e-1}^{r+e,r+1}=0$ we obtain that $E_e^{r,r}=\Ker(E_{e-1}^{r,r}\rightarrow E_{e-1}^{r-e,r-1})$.\\
If $r-p=e-1$ consider the maps
$$E_l^{p+l+1,r+1}\rightarrow E_l^{p,r}\rightarrow E_l^{p-l-1,r-1}$$
with $l\geq 1$.
We have that $E_l^{p+l+1,r+1}=0$ for all $1\leq l\neq e-1$ and that $E_l^{p-l-1,r-1}=0$ for all $l\geq 1$, thus $E_{e-1}^{p,r}=E_1^{p,r}$ and $E_{\infty}^{p,r}=E_e^{p,r}$.
Considering the maps $E_{e-1}^{p+e,r+1}\rightarrow E_{e-1}^{p,r}\rightarrow E_{e-1}^{p-e,r-1}$ we have that $E_{e-1}^{p-e,r-1}=0$ and so $E_e^{p,r}=E_{e-1}^{p,r}/ \Ima (E_{e-1}^{p+e,r+1}\rightarrow E_{e-1}^{p,r})$.\\
It follows that for all $m\geq 1$
$$E_\infty^{p,r} =
\begin{cases}
	0 & \text{if } (r-p\leq -1) \lor ((r-p\geq 1)\land(r-p\neq e-1)  )\\
	V_r & \text{if } r-p = 0\\
	W_{p,r} & \text{if } r-p=e-1
\end{cases},
$$
where $V_r$ is contained in  $E_1^{r,r}$ and $W_{p,r}$ is a quotient of $E_1^{p,r}$.\\
Consider now the filtration
$$H^i(\hat{X},\pi^*(mL)-t\hat{A})=F^0\supset F^1\supset \cdots \supset F^k\supset F^{k+1}=0$$
with $F^p/F^{p+1}=\Gr^p(H^i(\hat{X},\pi^*(mL)-t\hat{A}))=E_\infty^{p,i}$ for all $p\geq 0$.\\
Since $i> e-1$, then $F^p/F^{p+1}=0$ for all $p\neq i, i+1-e$. Thus the filtration is
$$H^i(\hat{X},\pi^*(mL)-t\hat{A})=F^0= F^{i+1-e}\supset F^{i+2-e}= F^i\supset F^{i+1}=0.$$
By the exact sequence
$$0\rightarrow F^{i+1} \rightarrow F^i\rightarrow F^i/F^{i+1}\rightarrow 0 $$
we observe that $F^{i+2-e}= F^i=V_i$. Moreover by the exact sequence
$$0\rightarrow F^{i+2-e} \rightarrow F^{i+1-e}\rightarrow F^{i+1-e}/F^{i+2-e}\rightarrow 0 $$
we have that 
$$h^i(\hat{X},\pi^*(mL)-t\hat{A})\leq \dim V_i + \dim W_{i+1-e,i}\leq h^i(X,mL-tA)+ \dim W_{i+1-e,i}.$$
Then we need to control the dimension of the $W_{i+1-e,i}$'s. To do this take $s\geq 0$ and consider the exact sequences
$$0\rightarrow \pi^*(mL-tA)+(t-s-1)E\rightarrow \pi^*(mL-tA)+(t-s)E \rightarrow (\pi^*(mL-tA)+(t-s)E)_{|E}\rightarrow 0$$
By projection formula
$$R^{e-1}\pi_*(\pi^*(mL-tA)+(t-s)E)=\mathcal{O}_X(mL-tA)\otimes R^{e-1}\pi_*((t-s)E) $$
while
$$R^{e-1}p_*((\pi^*(mL-tA)+(t-s)E)_{|E})=R^{e-1}p_*(i^*(\pi^*(mL-tA)+(t-s)E)))= $$
$$R^{e-1}p_*((p^*k^*(mL-tA)+(t-s)E_{|E}))= \mathcal{O}_Z((mL-tA)_{|Z})\otimes R^{e-1}p_*((t-s)E_{|E})).$$
Thus we get the exact sequence
\begin{equation}\label{merda}
0\rightarrow \mathcal{O}_X(mL-tA)\otimes R^{e-1}\pi_*((t-s-1)E)\rightarrow \mathcal{O}_X(mL-tA)\otimes R^{e-1}\pi_*((t-s)E)\rightarrow $$
$$ \rightarrow \mathcal{O}_Z((mL-tA)_{|Z})\otimes R^{e-1}p_*((t-s)E_{|E}))\rightarrow 0
\end{equation}
It follows that
$$\dim W_{i+1-e,i} \leq \dim E_1^{i+1-e,i}=h^{i+1-e}(X,\mathcal{O}_X(mL-tA)\otimes R^{e-1}\pi_*(tE))\leq$$
$$\leq\sum_{s=0}^{t-1} h^{i+1-e}(Z,\mathcal{O}_Z((mL-tA)_{|Z})\otimes R^{e-1}p_*((t-s)E_{|E})).$$
By \cite[Appendix A]{LazPAG} we have that
$$R^{e-1}p_*((t-s)E_{|E}) =
\begin{cases}
	0 & \text{if } s\geq t-e+1\\
	\Sym^{t-s-e}N_{Z/X} \otimes \det N_{Z/X} & \text{if } s\leq t-e
\end{cases}.
$$
It follows that $\dim W_{i+1-e,i}=0$ if $t\leq e-1$, while if $t\geq e$
$$\dim W_{i+1-e,i} \leq \sum_{s=0}^{t-1} h^{i+1-e}(Z,\mathcal{O}_Z((mL-tA)_{|Z})\otimes R^{e-1}p_*((t-s)E_{|E}))\leq$$
$$\leq \sum_{s=0}^{t-e} h^{i+1-e}(Z,\mathcal{O}_Z((mL-tA)_{|Z})\otimes \Sym^{t-s-e}N_{Z/X} \otimes \det N_{Z/X}).$$
Since $L$ is $q$-ample, there exists an integer $m_L>0$ such that
$$h^i(X,mL-tA)=0$$
for all $i>q$, $1\leq t\leq C$ and $m\geq m_L$.
Moreover, by Proposition \ref{credirr}, $L_{|Z}$ is also $q$-ample.
Hence for all $ e\leq t \leq C$ there exists an integer $m_t>0$ such that
$$h^{i+1-e}(Z,\mathcal{O}_Z((mL-tA)_{|Z})\otimes \Sym^{t-s-e}N_{Z/X} \otimes \det N_{Z/X})=0$$ 
for all $0\leq s\leq t-e$, $i>q+e-1$ and $m\geq m_t$.
It follows that
$$\dim W_{i+1-e,i} \leq \sum_{s=0}^{t-e} h^{i+1-e}(Z,\mathcal{O}_Z((mL-tA)_{|Z})\otimes \Sym^{t-s-e}N_{Z/X} \otimes \det N_{Z/X})=0 $$
for all $ e\leq t \leq C$, $i>q+e-1$ and $m>m_t$.\\
Then taking $m_{\pi^*L}=\max\{m_e,\dots ,m_{C},m_L\}$ we have the thesis.
\endproof

\begin{remark}
We remark that the previous result is sharp. Namely, we can not expect more reularity on $\pi^*L$. Indeed if $e=n-1$, then $\pi$ is an isomorphism and so if $L$ is $q$-ample but not $(q-1)$-ample, then $\pi^*L$ is $q$-ample but not $(q-1)$-ample.
\end{remark}

\end{document}